\documentclass[12pt]{article}
\usepackage{amssymb,amsthm,amsmath,latexsym}
\usepackage{url}
\usepackage{multirow}
\newtheorem{thm}{Theorem}

\newtheorem{lem}[thm]{Lemma}

\theoremstyle{remark}

\newcommand{\F}{\mathbb{F}}

\newcommand{\cC}{\mathcal{C}}
\newcommand{\cM}{\mathcal{M}}

\DeclareMathOperator{\Aut}{Aut}
\DeclareMathOperator{\PAut}{PAut}
\DeclareMathOperator{\GAut}{GAut}
\DeclareMathOperator{\MAut}{MAut}
\DeclareMathOperator{\wt}{wt}

\DeclareMathOperator{\supp}{supp}
\DeclareMathOperator{\gen}{gen}

\begin{document}

\title{Classification of binary self-dual $[76,38,14]$ codes with an automorphism of order $9$}
\date{}
\author{
Nikolay Yankov,
Radka Russeva and Emine Karatash
\thanks{
Faculty of Mathematics and Informatics,
Konstantin Preslavski University of Shumen,
Shumen, 9712, Bulgaria.}
}

\maketitle

\begin{abstract}
\small{
Using the method for constructing binary self-dual codes with an automorphism of order square of a prime number we
have classified all binary self-dual codes with length $76$ having minimum distance $d=14$ and automorphism of order $9.$
Up to equivalence, there are six self-dual $[76,38,14]$ codes with an automorphism of type $9$-$(8,0,4).$
All codes obtained have new values of the parameter in their weight enumerator thus more than doubling the number of known values.

}
\end{abstract}


\section{Introduction}

In this paper, we are interested in the classification of the extremal binary self-dual $[76,38,14]$ codes with an
automorphism of order $9.$ It was motivated by the following reasons.

Firstly, there are only three known extremal binary self-dual $[76,38,14]$ codes, constructed by Dontcheva and Yorgov via an automorphism of order 19 \cite{Dontcheva2003a}.
These three codes are not only shadow optimal but also shortest known self-dual code with minimal distance 14.
One of these three codes was the first ever found in the literature and it was discovered by Baartmans and Yorgov \cite{Baartmans2003}.

Secondly, Bouyuklieva, et al \cite{Bouyuklieva2005}
presented a method for constructing binary self-dual codes with an automorphism of order $p^2$ and classified all optimal binary self-dual
codes self-dual codes of lengths $44\leq n\leq 54$ having an automorphism of order $9$. The case for the length of an optimal
binary self-dual code with automorphism of such order was considered by Yankov in \cite{Yankov2012} where it was proved that a doubly-even self-dual
$[72, 36, 16]$ codes with an automorphism of order 9 does not exists.

A linear $[n, k]$ \emph{code} $C$ is a $k$-dimensional subspace of the vector space $\F_q$,
where $\F_q$ is the finite field of $q$ elements. The elements of $C$ are called \emph{codewords},
and the \emph{(Hamming) weight} of a codeword $v\in C$ is the number of the non-zero coordinates of $v$.
We use $\wt(v)$ to denote the weight of a codeword.
The \emph{minimum weight} $d$ of $C$ is the smallest weight among all its non-zero codewords,
and $C$ is called an $[n, k, d]_q$ code. A matrix whose rows form a basis of $C$ is called
a \emph{generator matrix} of this code and we denote this by $\gen(C).$ Every code satisfies the Singleton bound $d\leq n-k+1$.
A code is \emph{maximum distance separable} or \emph{MDS} if $d=n-k+1$, and \emph{near MDS} or \emph{NMDS} if $d=n-k$.

For every $u = (u_1,\dots, u_n)$ and $v = (v_1,\dots, v_n)$ from $\F_2^n$, $u.v=\sum\limits_{i=1}^{n}{u_iv_i}$
defines the \emph{inner product} in $\F_2^n$. The \emph{dual code} of $C$ is $C^\perp=\{v\in \F_2^n\;\vert\; u.v = 0, \forall\; u\in C\}$.
If $C\subset C^\perp$, $C$ is called \emph{self-orthogonal}, and if $C=C^\perp$, we say that $C$ is \emph{self-dual}.
We call a binary code \emph{self-complementary} if it contains the all-ones vector. Every binary self-dual code is self
complementary.

A self-dual code is \emph{doubly even} if all codewords have weight divisible by four, and \emph{singly even} if there is at least
one nonzero codeword $v$ of weight $\wt(v)\equiv 2(\text{mod}\;4)$. Self-dual doubly even codes exist only if $n$ is a multiple of eight.

The \emph{Hermitian inner product} on $\F_4^n$ is given by $u.v=\sum\limits_{i=1}^{n}{u_iv_i^2}$ and we denote by
$C^{\perp H}$ the dual of $C$ under Hermitian inner product. $C$ is \emph{Hermitian self-dual} if $C=C^{\perp H}$.

The weight enumerator $W(y)$ of a code $C$ is defined as $W(y)=\sum_{i=0}^n A_iy^i$, where $A_i$ is the number
of codewords of weight $i$ in $C$. Following \cite{Huffman2003} we say that two linear codes $C$ and $C'$ are
\emph{permutation equivalent} if there is a permutation of coordinates which sends $C$ to $C'$.
The set of coordinate permutations that maps a code $C$ to itself forms a group denoted by $\PAut(C)$.
Two codes $C$ and $C'$ of the same length over $\F_q$ are \emph{equivalent} provided there is a monomial
matrix $M$ and an automorphism $\gamma$ of the field such that $C=C'M\gamma$.
The field $\F_4$ has an automorphism $\gamma$ given by $\gamma(x)=x^2$.

The set of monomial matrices that maps $C$ to itself forms the group $\MAut(C)$ called the
\emph{monomial automorphism group} of $C$. The set of maps of the form $M\gamma$, where $M$ is a monomial matrix and
$\gamma$ is a field automorphism, that sends $C$ to itself, forms the group $\GAut(C)$, called the automorphism group of $C$.
In the binary case all three groups are identical. In general, $\PAut(C)\subseteq \MAut(C) \subseteq \GAut(C)$.

An automorphism $\sigma\in{\cal S}_n$, $|\sigma|=p^2$ is of \emph{type} $p^2$-$(c,t,f)$
if when decomposed to independent cycles it has $c$ cycles of length $p^2$,
$t$ cycles of length $p$, and $f$ fixed points. Obviously, $n=cp^2+tp+f$.

This paper is organized in the following way. First in Section \ref{constr} we introduce to the reader the main results about the method we use.
Section \ref{sec76} shows the application of the method and the construction of six new binary self-dual $[76, 38, 14]$ codes.

\section{Construction Method}\label{constr}

In \cite{Bouyuklieva2005} a method for constructing binary self-dual codes having an automorphism of order $p^2$, where $p$
is an odd prime, was presented. We consider the case $p=3.$

Let $C$ be a self-dual $[76, 38, 14]$ code having an automorphism $\sigma$ of type $9$-$(c,t,f).$ 
In \cite{Bouyuklieva2007a} (Lemma 6) it is proved that $\sigma$ is of type $9$-$(8,0,4)$, i.e. $c=8$, $t=0$ and $f=4$.
Thus we have
\begin{equation}\label{sigma1}
\sigma=(1,2,\dots,9)(10,11,\dots,18)\dots(64,65,\dots,72)(73)\dots(76).
\end{equation}
Denote by $\Omega_i$, $i=1,\dots,12$ the cycles in $\sigma$.
Define $${F_\sigma }(C)=\{ v \in C\;\vert\;v\sigma  = v\},$$
$${E_\sigma }(C) = \{ v \in C\;\vert\;\wt(v\vert{\Omega _i}) \equiv 0\pmod{2}\},$$
where $v\vert\Omega _i$ denotes the restriction of $v$ to $\Omega_i$.
Clearly $v\in F_\sigma(C)$ iff $v\in C$ is constant on each cycle. Denote
$\pi : F_\sigma(C)\to \F_2^{12}$ the projection map where if $v\in F_\sigma(C),$
$(\pi(v))_i=v_j$ for some $j\in\Omega_i$, $i=1,\dots, 12$.
Then the following lemma holds.
\begin{lem}\cite{Bouyuklieva2005} $C=F_\sigma(C)\oplus E_\sigma(C)$. $C_\pi=\pi(F_\sigma(C))$ is a
binary self-dual code of length $12.$
\end{lem}

Thus each choice of the codes $F_\sigma(C)$ and $E_\sigma(C)$ determines a self-dual code $C$. So for a given length all
self-dual codes with an automorphism $\sigma$ can be obtained.

Denote with $E_\sigma(C)^{\ast}$ the subcode $E_\sigma(C)$ with the last 4 zero coordinates deleted. $E_\sigma(C)^{\ast}$ is a self-orthogonal binary code of length $8.3^2=72$ and dimension $\frac{8}{2}(3^2-1)=32$.
For $v\in {E_\sigma (C)}^{\ast}$ we let $v\vert\Omega_i=(v_0,v_1,\cdots,v_{8})$ correspond to the
polynomial $v_0+v_1 x+\dots+v_{8}x^{8}$ from ${\cal T}$, where ${\cal T}$ is the ring of even-weight polynomials in $\F_2
[x]/(x^9-1)$. Thus we obtain the map $\varphi:E_\sigma(C)^{\ast}\to {\cal T}^8 $. Denote $C_\varphi=\varphi(E_\sigma(C)^{\ast})$.

Let $e_1=x^8+x^7+x^5+x^4+x^2+x$ and $e_2=x^6+x^3$. In our work \cite{Bouyuklieva2005} we proved that
${\cal T}=I_1\oplus I_2$, where $I_1=\{0,e_1,\omega=xe_1,\overline{\omega}=x^2e_1\}$ is a field with identity $e_1$
and $I_2$ is a field with $2^6$ elements with identity $e_2$. The element $\alpha=(x+1)e_2$ is
a primitive element in $I_2$ so $I_2=\{0,\alpha^k, 0\leq k\leq 62\}$.

The following theorem is from \cite{Bouyuklieva2007a}.

\begin{thm}{\cite{Bouyuklieva2007a}}\label{Theorem1} $C_\varphi=M_1\oplus M_2$, where
$M_j=\{u\in E_\sigma(C)^{\ast}|u_i\in I_j, i=1,\dots,8\}$, $j=1,2$. Moreover
$M_1$ and $M_2$ are Hermitian self-dual codes over the fields $I_1$ and $I_2$, respectively.
If $C$ is a binary self-dual code having an automorphism $\sigma$ of type
\eqref{sigma1} then $E_\sigma(C)^{\ast}=E_1\oplus E_2$ where $M_i=\varphi(E_i)$, $i=1,2$.
\end{thm}

This proves that $C$ has a generator matrix of the form
\begin{equation}\label{genmat}
{\cal G}=
\left(\begin{array}{c}
\varphi^{-1}(M_2) \ 0 \ 0 \ 0 \ 0 \\
\varphi^{-1}(M_1) \ 0 \ 0 \ 0 \ 0 \\
F_\sigma\\
\end{array}\right).
\end{equation}

Let $B_s$ and $E_s$ denote the number of words of weight $s$ is $F_\sigma(C)$ and $E_\sigma(C)^{\ast}$, respectively.
Every word of weight $s$ in $E_\sigma(C)^{\ast}$ is in an orbit of length $3,$ therefore, $E_s\equiv 0\pmod{3}$ and $A_s\equiv B_s\pmod{3}$ for $1\leq s\leq n.$

Since the minimum distance of $C$ is 16 the code $M_2$ is a $[8, 4]$ Hermitian self-dual code over $\F_{64}$, having minimal distance $d\geq 4$. Using Singleton bound $d\leq n-k+1$ we have $d=5$  or  $d=4$.
The case $d=5$ is studied in \cite{Bouyuklieva2007a} and there are exactly 96 MDS Hermitian $[8, 4, 5]_{64}$ self-dual codes
such that the minimum distance of $\varphi^{-1}(M_2)$ is 16. The case for the near MDS codes is completed in \cite{Yankov2012} and the number of the codes is 26 and we state the following.

\begin{thm}[\cite{Bouyuklieva2007a}, \cite{Yankov2012}]
Up to equivalence, there are exactly $122$ Hermitian $[8, 4]_{64}$ self-dual codes such that the minimum distance of $\varphi^{-1}(M_2)$ is $16.$
\end{thm}

We denote these codes by $\cM_{2,i}$ for $1\leq i\leq 122.$ Their generator parameters can be obtained from \cite{Yankov2012}.

We fix the upper part of $\cal{G}$ in \eqref{genmat} to be generated by one of the 122 already constructed
Hermitian MDS or NMDS $[8,4]$ codes. Now we continue with construction of the middle part,
i.e. the code $M_1$. Theorem \ref{Theorem1} states that
$M_1$ is a quaternary Hermitian self-dual $[8, 4, 4]$ code. There exists
a unique such code $e_8$ \cite{Conway1979} with a generator matrix
$Q_1=\left(\begin{array}{c}
10000111 \\
01001011 \\
00101101 \\
00011110 \\
\end{array}\right)$. We have to put together the two codes from $M_2$ and $M_1$ in \eqref{genmat}, but
we have to examine carefully all transformations on $Q_1$ that can lead to a different joined code.
The full automorphism group of $e_8$ is of order $2.3^8(8!)$
and we have to consider the following transformations that preserve the decomposition of the code $C:$
\begin{itemize}
  \item[\em{(i)}] a permutation $\tau\in S_8$ acting on the set of columns.
  \item[\em{(ii)}] a multiplication of each column by a nonzero element $e_1,\omega$ or $\overline{\omega}$ in $I_1$.
  \item[\em{(iii)}] a Galois automorphism $\gamma$ which interchanges $\omega$ and $\overline{\omega}$.
\end{itemize}

The action of (i) and (ii) can be represented by a monomial matrix $M=PD$ for a
diagonal matrix $D$ and permutational matrix $P$. Since every column of $Q_1$
consists only of $0$ and $1$ the action of $PD\gamma$ on $Q_1$ can be obtained via
$P\overline{D}$. Thus we apply only transformations (i) and (ii).

Denote by $M_1^{\tau}$ the code determined by the matrix $Q_1$ with columns permuted by
$\tau$. To narrow down the computations we can use $\PAut({M_1})=\langle\,$(47)(56), (45)(67),
(12)(3586), (24)(68), (34)(78)$\,\rangle,$ $|\PAut({M_1})|=1344$ and the right transversal $T$ of
$S_8$ with respect to $\PAut({M_1})$
\begin{align*}
T{=}&
\{(), (78), (67), (678), (687), (68), (56), (56)(78), (567), (5678), (5687), \\
&(568), (576), (5786), (57), (578), (57)(68), (5768), (5876), (586), (587), (58), \\
&(5867), (58)(67), (45678), (4568), (4578), (45768), (458), (458)(67)\}.
\end{align*}

For every one of the 122 codes $\cM_{2,i}$ and $\tau\in T$ we considered $3^8$ possibilities for $\gen(M_1^\tau)$
and checked the minimum distance in the corresponding binary code $E_\sigma(C)^{\ast}$. We state the following result.

\begin{thm}\label{Thm_72}
There are exactly $36659$ inequivalent self-orthogonal $[72, 32, 16]$ codes having an automorphism with $8$ cycles of order $9.$
\end{thm}

Denote the codes obtained by $C_{72,i}$, $i=1,\dots,36659.$  In Table \ref{n72_aut} and Table \ref{n72_wt} we summarize the values of the order of the automorphism groups  $|\Aut|$ and the number $A_{16}$ of codewords of weight 16 for these codes.

\begin{table}[htb]
\centering
\caption{The cardinality of the automorphism groups of the $[72, 32, 16]$ codes}\label{n72_aut}
\begin{tabular}{|c|c|c|c|c|c|c|}
\hline
$|\Aut|$&9&18&27&36&54&72\\
\hline
\# of codes&35876&730&24&25&2&2\\
\hline
\end{tabular}
\end{table}

\begin{table}[htb]
\centering
\caption{Number of inequivalent $[72, 32, 16]$ codes with $A_{16}$}\label{n72_wt}
{\footnotesize\begin{tabular}{|c|c||c|c||c|c||c|c||c|c||c|c|}
\hline
$A_{16}$&\#&$A_{16}$&\#&$A_{16}$&\#&$A_{16}$&\#&$A_{16}$&\#&$A_{16}$&\#\\
\hline
\hline
14751&1&14967&454&15183&806&15399&207&15615&27&15831&2\\
\hline
14760&2&14976&479&15192&787&15408&212&15624&28&15840&1\\
\hline
14769&1&14985&569&15201&740&15417&180&15633&26&15849&1\\
\hline
14778&2&14994&598&15210&798&15426&193&15642&20&15858&1\\
\hline
14787&6&15003&635&15219&654&15435&148&15651&20&15867&1\\
\hline
14796&8&15012&722&15228&674&15444&161&15660&8&15876&1\\
\hline
14805&12&15021&740&15237&654&15453&145&15669&9&15894&1\\
\hline
14814&17&15030&760&15246&687&15462&118&15678&8&15903&1\\
\hline
14823&25&15039&764&15255&615&15471&127&15687&15&15912&3\\
\hline
14832&32&15048&787&15264&521&15480&120&15696&6&15921&1\\
\hline
14841&45&15057&807&15273&544&15489&116&15705&13&15948&1\\
\hline
14850&62&15066&826&15282&503&15498&102&15714&9&15957&1\\
\hline
14859&70&15075&815&15291&504&15507&75&15723&8&15966&2\\
\hline
14868&93&15084&889&15300&424&15516&75&15732&8&15975&1\\
\hline
14877&127&15093&910&15309&446&15525&68&15741&7&15984&1\\
\hline
14886&155&15102&860&15318&428&15534&56&15750&5&16011&1\\
\hline
14895&179&15111&962&15327&416&15543&60&15759&8&16029&1\\
\hline
14904&213&15120&832&15336&385&15552&48&15768&9&16200&1\\
\hline
14913&290&15129&827&15345&357&15561&48&15777&1&16218&1\\
\hline
14922&264&15138&863&15354&340&15570&39&15786&2&16632&1\\
\hline
14931&317&15147&862&15363&345&15579&48&15795&5&16848&1\\
\hline
14940&326&15156&855&15372&288&15588&38&15804&3&17604&1\\
\hline
14949&401&15165&847&15381&267&15597&39&15813&2&&\\
\hline
14958&419&15174&784&15390&233&15606&30&15822&3&&\\
\hline
\end{tabular}}
\end{table}

\section{Construction of new $[76,38,14]$ codes with an automorphism of type $9$-$(8,0,4)$}\label{sec76}

The highest attainable minimum weight for length 76 is 14 and there are three possible weight enumerators and shadows \cite{Dougherty1997}:
\begin{align*}
&\left\{\begin{array}{r@{\,}l}
W_{76,1}&=1+(4750-16\alpha)y^{14}+(79895+64\alpha)y^{16}+(915800+ 64\alpha)y^{18}+\cdots\\
S_{76,1}&=\alpha y^{10}+(9500-14\alpha)y^{14}+(1831600+91\alpha)y^{18}+\cdots\\
        &(0\leq\alpha\leq 296)\\
\end{array}
\right.\\
&\left\{\begin{array}{r@{\,}l}
W_{76,2}&=1+2590y^{14}+106967y^{16}+674584y^{18}+\cdots\\
S_{76,2}&=y^2+8954y^{14}+1836865y^{18}+105664452y^{22}+\cdots\\
\end{array}
\right.\\
&\left\{\begin{array}{r@{\,}l}
W_{76,3}&=1+(4750+16\alpha)y^{14}+(80919-64\alpha)y^{16}+(905560-64\alpha)y^{18}+\cdots\\
S_{76,3}&=y^6+(-16-\alpha)y^{10}+(9620+14\alpha)y^{14}+(1831040-91\alpha)y^{18}+\cdots\\
&(-296\leq\alpha\leq -16)\\
\end{array}
\right.
\end{align*}
There are only three known codes with $\alpha=0$ for $W_{76,1}$ \cite{Dontcheva2003a}, possessing an automorphism of type $19$-$(4,0).$

Now $C_\pi$ is a binary self-dual $[12,6]$ code. Up to equivalence there are three such codes $6i_2, 2i_2+h_8$ and $d_{12}$ \cite{Pless1972}.
In the case of $6i_2$ we can not fix any point since then there will be a codeword of weight 10 in $C$.
When $C_\pi\cong 2i_2+h_8$ we have to take the four fixed points from the $h_8$ summand. Since the automorphism group of $h_8$ is 3-transitive
we can take any three points from it

and we have to choose one more cyclic point from the last five. We checked all five different splits and found a vector in $F_\sigma(C)$ with weight $d<14.$
Lastly, when $C_\pi\cong d_{12}$, for every 4-weight codeword we have to choose at least two coordinates from its support.

The code $d_{12}$ possesses a cluster $\{\{1,2\}, \{3,4\},\{5,6\}, \{7,8\}, \{9,10\}, \{11,12\}\}$ so we have to choose the four fixed points from different duads.
Up to a permutation of the cyclic points or a permutation of the fixed points we have a unique generating matrix
$$G_2=\left(\begin{array}{cccccccc|cccc}
1&1&0&0&0&0&0&0&1&1&0&0\\
0&1&1&0&0&0&0&0&0&1&1&0\\
0&0&1&1&0&0&0&0&0&0&1&1\\
0&0&0&1&1&1&0&0&0&0&0&1\\
0&0&0&0&1&1&1&1&0&0&0&0\\
1&1&1&1&0&1&0&1&0&0&0&0
\end{array}\right)$$ for the code $C_\pi.$

By Q-extensions \cite{Bouyukliev2007a} we obtained  $G''=\left\langle (1,2), (2,4,3)(5,7)(6,8), (5,6)(7,8) \right\rangle$ the subgroup of the symmetric group $S_8$
that preserves the code generated by $G_2.$ The group $G''$ has cardinality 420. To construct a generator matrix of a $[76,38]$ self-dual code in form \eqref{genmat} we fix a generator matrix of $E_\sigma(C)^{\ast}$ and we use the matrix  $G_2$ with columns permuted by $\mu$ for all permutations $\mu \in G''$.

Our exhaustive search gives the following result.

\begin{thm}
\label{n76}
Up to equivalence there exist exactly $6$ binary self-dual $[76, 38, 14]$ codes with an automorphism of type $9$-$(8,0,4).$
All codes have weight enumerators $W_{76,1}$ for $\alpha=4$ or $13$ and automorphism groups of order 9.
\end{thm}

The generator parameters and the weight enumerator for the six binary self-dual $[76, 38, 14]$ codes, denoted by  $\cC_{76,i}$ $1\leq i\leq 6$, are displayed in Table \ref{n76t}.
The notation $\tau, D$ in Table \ref{n76t} means that we are using the permutation $\tau\in T$ on $M_1^\tau$ and then a multiplication of each column by the corresponding element in $D.$
Alternatively the generator matrices of the codes $\cC_{76,i}$ for $i=1,2,\dots,6$ can be obtained online at ``http://shu.bg/tadmin/upload/storage/2599.txt''.

\begin{table}[htb]
\centering
\caption{The generators for the constructed $[76, 38, 14]$ codes}\label{n76t}
{\footnotesize\begin{tabular}{|c|c|c|c|c|}
\hline
code&$\cC_{72,i}$&$\tau,D$&$\supp(C_\pi)$&$\alpha$\\
\hline
\hline
\multirow{2}{*}{$\cC_{76,1}$}&\multirow{2}{*}{11}&$(4,5,7,8),$&$\{1,3,9,10\},\{3,5,10,11\},\{5,8,11,12\},$&\multirow{2}{*}{$4$}\\
&&$(1,1,\overline{\omega},\overline{\omega},{\omega},1,1,{\omega})$&$\{2,7,8,12\},\{2,4,6,7\},\{1,3,5,6,7,8\}$&\\
\hline
\multirow{2}{*}{$\cC_{76,2}$}&\multirow{2}{*}{11}&$(4,5,7,8),$&$\{1,3,9,10\},\{3,5,10,11\},\{5,8,11,12\},$&\multirow{2}{*}{$4$}\\
&&$(1,1,\overline{\omega},\overline{\omega},{\omega},1,1,{\omega})$&$\{2,7,8,12\},\{2,4,6,7\},\{1,3,4,5,7,8\}$&\\
\hline
\multirow{2}{*}{$\cC_{76,3}$}&\multirow{2}{*}{36}&$(4,5,8),$&$\{2,4,9,10\},\{4,6,10,11\},\{6,8,11,12\},$&\multirow{2}{*}{$4$}\\
&&$(1,\overline{\omega},\overline{\omega},\overline{\omega},{\omega},\overline{\omega},1,\overline{\omega})$&$\{1,7,8,12\},\{1,3,5,7\},\{2,4,5,6,7,8\}$&\\
\hline
\multirow{2}{*}{$\cC_{76,4}$}&\multirow{2}{*}{36}&$(4,5,8),$&$\{2,4,9,10\},\{4,6,10,11\},\{6,8,11,12\},$&\multirow{2}{*}{$4$}\\
&&$(1,\overline{\omega},\overline{\omega},\overline{\omega},{\omega},\overline{\omega},1,\overline{\omega})$&$\{1,7,8,12\},\{1,3,5,7\},\{2,3,4,6,7,8\}$&\\
\hline
\multirow{2}{*}{$\cC_{76,5}$}&\multirow{2}{*}{106}&$(4,5,6,7,8),$&$\{3,4,9,10\},\{4,5,10,11\},\{5,7,11,12\},$&\multirow{2}{*}{$13$}\\
&&$(1,{\omega},1,{\omega},{\omega},{\omega},{\omega},1)$&$\{1,6,7,12\},\{1,2,6,8\},\{3,4,5,6,7,8\}$&\\
\hline
\multirow{2}{*}{$\cC_{76,6}$}&\multirow{2}{*}{106}&$(4,5,6,7,8),$&$\{3,4,9,10\},\{4,5,10,11\},\{5,7,11,12\},$&\multirow{2}{*}{$13$}\\
&&$(1,{\omega},1,{\omega},{\omega},{\omega},{\omega},1)$&$\{1,6,7,12\},\{1,2,6,8\},\{2,3,4,5,6,7\}$&\\
\hline
\end{tabular}}
\end{table}

\section*{Acknowledgements}

This work was supported by Shumen University under Project No RD-08-68/02.02.2017.

\end{document}